\documentclass[11pt]{article}
\usepackage{amsthm,amsmath,amssymb,amsfonts}
\newcommand{\eqdef}{\stackrel{{\rm}}{=}}

\newcommand{\cw}{\stackrel{D}{\rightharpoonup}}

\newcommand{\N}{\mathbb{N}}
\newcommand{\R}{\mathbb{R}}

\newcommand{\Z}{\mathbb{Z}}
\newcommand{\wlim}{\operatorname{w-lim}}

\begin{document}

% theorem style plan --default
\newtheorem{theorem}{THEOREM}[section]
\newtheorem{corollary}[theorem]{COROLLARY}
\newtheorem{lemma}[theorem]{LEMMA}
\newtheorem{proposition}[theorem]{PROPOSITION}

\theoremstyle{definition}
\newtheorem{definition}[theorem]{Definition}
\newtheorem{conjecture}{Conjecture}
\newtheorem{example}{Example}
\newtheorem{problem}{Problem}

\theoremstyle{remark}
\newtheorem{remark}[theorem]{Remark}
\newtheorem{remarks}[theorem]{Remarks}
\newtheorem{claim}[theorem]{Claim}  
\newtheorem{summary}[theorem]{Summary}   \renewcommand{\thesummary}{}
\newtheorem{case}[theorem]{Case}
\newtheorem{ack}[theorem]{Acknowledgment}  \renewcommand{\theack}{}

\begin{titlepage}
\title{
{\bf Cocompact imbeddings and structure of weakly convergent sequences}
}
\author{
\\K.Tintarev
\\Uppsala University,
\\tintarev@math.uu.se\\
}
\date{}\maketitle
\end{titlepage} 
\begin{abstract} 
Concentration compactness method is a powerful techniques for establishing existence of minimizers for inequalities and of critical points of functionals in general. The paper gives a functional-analytic formulation for the method in Banach space,
generalizing the Hilbert space case elaborated in \cite{ccbook}. The key object is a dislocation space - a triple $(X,F,D)$, where $F$ is a convex functional that defines a norm on  Banach space $X$, and $D$ is a group of isometries on $X$. Bounded sequences in dislocation spaces admit a decomposition into an asymptotic sum ``profiles" $w^{(n)}\in X$ dislocated by actions of $D$, that is, a sum of the form 
$\sum_ng^{(n)}_kw^{(n)}$, $g^{(n)}_k\in D$, while the remainder term converges weakly under actions of any sequence $g_k\in D$ ({\em $D$-weak convergence}). This decomposition allows to extend the weak convergence argument from variational problems with compactness to problems where $X$ is {\em cocompactly} (relatively to the group $D$) imbedded into a Banach space $Y$, that is, when every sequence $D$-weakly convergent in $X$ is convergent in the norm of $Y$. We prove a general statement on existence of minimizers in cocompact imbeddings that applies, in particular to Sobolev imbeddings which lack compactness (unbounded domain, critical exponent) including the subelliptic Sobolev spaces and spaces over Riemannian manifolds.
\end{abstract}
\section{Introduction}
Minimizers for an inequality in a functional space often do not exist, or cannot be obtained by a straightforward compactness reasoning, since in general one may expect that the minimization sequence converges only weakly. Concentration compactness method presented in the celebrated papers of P.-L.Lions \cite{PLL1a, PLL1b, PLL2a, PLL2b}) uses detailed structural information about minimization sequences in order to verify convergence in problems that a priori lack compactness. The core idea of concentration compactness is that if the problem possesses a noncompact invariance group $G$, lack of convergence can be attributed to the action of $G$, and thus a given sequence becomes convergent only after the terms (``profiles"), dislocated by the transformations, are ``factored out". Elaborations of the original classification of weak convergent sequences by Lions into tight, vanishing and dichotomous, which are often called ``splitting lemmas" were given for specific cases by Struwe \cite{Struwe84}, Brezis and Coron \cite{BrezisCoron}, Lions himself \cite{Lions87}, and numerous authors afterwards. The ``splitting lemmas", which were originally established for critical sequences of specific functionals in specific functional spaces, have been later summarized by the author in a structural statement that holds in the general Hilbert space (see \cite{ccbook} and references therein) using asymptotic orthogonality of dislocated profiles: if $g_k\in D$, where $D$ is a fixed group of unitary operators, $u_k\rightharpoonup w$,  $g_ku_k\rightharpoonup w_2$, and $v_k\eqdef u_k-w_1-g_k^{-1}w_2$ then $u_k=w_1+g_k^{-1}w_2+v_k$ is an asymptotically orthogonal sum in the sense that the scalar product of any two terms in it converges to zero. Furthermore,
this construction may be iterated. Under general conditions, subtraction of all dislocated profiles of a bounded sequence (nonzero weak limits of sequences $g_ku_k$ with different sequences $g_k$) amounts to a sequence that weakly converges to zero under all dislocations ($D$-weak convergence).
In fact this construction is useful only to an extent that $D$-weak convergence
is meaningful. One may say that the Hilbert space is cocompactly (relatively to $D$) imbedded into a Banach space $Y$ if $D$-weak convergence in $X$ implies convergence in $Y$. For example, subcritical Sobolev imbeddings on complete Riemannian manifolds are cocompact with respect to the action of any subgroup of the isometry group of the manifold if the manifold itself is cocompact with respect to this subgroup.

In the present paper we give a tentative generalization of this framework to Banach spaces, where one can no longer rely on the notion of asymptotic orthogonality. Its natural counterpart is asymptotic additivity or subadditivity of energy functionals with respect to dislocated profiles (it makes sense indeed to call a functional with such additivity property an energy, indicating that it is asymptotically additive over asymptotically separate (e.g. with asymptotically disjoint supports) clusters of the physical system that it models). Such asymptotic additivity is realized, in particular, in Brezis-Lieb lemma \cite{BL} 

Many applications of the concentration compactness method, such as existence of minimizers in isoperimetric problems or compactness of imbeddings of subspaces of functions with symmetries, are realized already on the functional-analytic level, with immediate applications to Sobolev spaces $W^{m,p}$ over Riemannian (and sub-Riemannian) manifolds and their flask subdomains. 

In Section 2 we prove the main structural theorem. Section 3 deals with functional-analytic statements on existence of minimizers in isoperimetric problems. Section 4 extends the results of two previous cases to non-invariant subspaces, and in Section 5 some compactness results are given.

\section{Dislocation space and weak convergence decomposition}
In this section we prove a structural theorem for bounded sequences in a class of Banach spaces associated with convex functionals. 
\begin{lemma}
\label{prop:OS}
Let $X$ be a vector space and let $F\in C^2(X)$
be an even nonnegative convex function with $F^{-1}(0)=\{0\}$. 
Then the map $\lambda:X\to[0,\infty)$,
\begin{equation}
\label{Orl-norm}
\lambda(u)=\inf\{\lambda>0: F(\lambda^{-1}u)\le 1\}
\end{equation}
is a norm on $X$  and
$\lambda=\|u\|$, for any $u\in X\setminus\{0\}$, is the unique solution of $F(\lambda^{-1}u)=1$.
%Furthermore, the function $F$ is continuous at zero.
\end{lemma}
\begin{proof} 
Homogeneity of $\lambda(u)$ is immediate from the definition.
If $u=0$ then $F(\lambda^{-1}u)=0$ for all $\lambda>0$ and thus $\lambda(0)=0$.
Since $F^{-1}(0)=\{0\}$, for every $u\in X\setminus\{0\}$, the even convex function $t\in\R\mapsto F(tu)$ is strictly monotone and unbounded from above.
In particular, $\lambda(u)>0$ whenever $u\neq 0$. Furthermore, by strict monotonicity, $F(\lambda^{-1}u)=1$ has a unique solution $\lambda_1$, and since $F(\lambda^{-1}u)>1$ whenever $\lambda<\lambda_1$, the infimum in \eqref{Orl-norm} is attained at $\lambda_1=\lambda(u)$.
It remains to prove the triangle inequality.
By convexity of $F$, 
\begin{eqnarray*}
F\left(\frac{u+v}{\lambda(u)+\lambda(v)}\right)=
F\left(\frac{\lambda(u)}{\lambda(u)+
\lambda(v)}\frac{u}{\lambda(u)}+
\frac{\lambda(v)}{\lambda(u)+\lambda(v)}\frac{v}{\lambda(v)}\right)
&\le&\\
\frac{\lambda(u)}{\lambda(u)+\lambda(v)}F\left(\frac{u}{\lambda(u)}\right)+
\frac{\lambda(v)}{\lambda(u)+\lambda(v)}F\left(\frac{v}{\lambda(v)}\right)&=&
\\\frac{\lambda(u)}{\lambda(u)+\lambda(v)}+
\frac{\lambda(v)}{\lambda(u)+\lambda(v)}&=&1.
\end{eqnarray*}
%Let now $u_k\in X$ be such that $\lambda(u_k)\to 0$. Then
%\begin{eqnarray*}
%F(u_k)=F\left(\lambda(u_k)\frac{u_k}{\lambda(u_k)}+(1-\lambda(u_k))0\right)&\le&\\ \lambda(u_k)F\left(\frac{u_k}{\lambda(u_k)}\right)+
%(1-\lambda(u_k))F(0)&= &\lambda(u_k)\to 0,
%\end{eqnarray*}
%which proves the continuity at zero.
\end{proof}

\begin{definition}
\label{DS} A dislocation space is a triple (X,F,D), where 
the pair $(X,F)$ is as in Lemma~\ref{prop:OS}, $F\in C(X)$ is uniformly continuous on bounded sets, the Banach space $X$ is separable and reflexive, and $D$ is a group of linear operators on $X$, closed with respect to the strong (elementwise) convergence, satisfying 
$F\circ g=F$ for all $g \in D$,  and such that
\begin{equation}
\label{newii}
g_k\in D, g_k\not\rightharpoonup 0, u_k\rightharpoonup 0 \Rightarrow g_k u_k \rightharpoonup 0 \text{ on a subsequence}.
\end{equation}
Moreover, if sequences $\{g_k^{(n)}\}_{k\in\N}\subset D$, $n=1,\dots,M$, $M\in\N$, satisfy 
\begin{equation}
\label{separates}
{g_k^{(m)}}^{-1}{g_k^{(n)}}\rightharpoonup 0,\quad m\neq n,
\end{equation}
and $u_k\in X$ is a bounded sequence such that
${g_k^{(n)}}^{-1}u_k\rightharpoonup w^{(n)}$, $n=1,\dots,M$, then
\begin{equation}
\label{wBL}
\liminf F(u_k)\ge\sum_{n=1}^MF(w^{(n)}).
\end{equation}
and 
\begin{equation}
\label{w1BL}
F(\sum_{n=1}^Mg_k^{(n)}w^{(n)})\to \sum_{n=1}^MF(w^{(n)}).
\end{equation}
\end{definition}
\begin{remark}
It is easy to see that conditions \eqref{wBL} and \eqref{w1BL} are satisfied if $F$ satisfies the Brezis-Lieb property: 
$$u_k\rightharpoonup u\Rightarrow F(u_k)-F(u)-F(u_k-u)\to 0.$$
In particular, if $X$ is a Hilbert space and $F(u)=\|u\|^2$, then
$$\|u_k\|^2-\|u_k-u\|^2-\|u\|^2=2(u_k,u)-2\|u\|^2\to 0.$$
When $F(u)=\int \varphi(u)d\mu$ with $\varphi$ from a class of functions on a measure space that includes $\varphi(t)=|t|^p$, $p\in(1,\infty)$, Brezis-Lieb property has been verified in \cite{BL} under an additional condition $u_k\to 0$ a.e., although since $L^2$ is a Hilbert space, this condition redundant when $p=2$.
\end{remark}
\noindent Examples of dislocation spaces:
\begin{enumerate}
\item $(H,F,D)$, where $H$ is a separable Hilbert space; $F(u)=\|u\|^2$; and the group $D$ of unitary operators satisfies \eqref{newii}, in particular, as in any of the examples below with $p=2$. This case is elaborated in \cite{ccbook}.
\item $(W^{1,p}(M),\|\cdot\|^p,D)$, where $M$ is a complete 
(~sub-~)~Riemannian manifold, $W^{1,p}(M)$, $p>1$, is a Sobolev space associated with the $p$-(~sub-~)Laplacian and $D=\{u\mapsto u\circ\eta\}_{\eta\in Iso(M)}$.
\par
In particular,  $(W^{1,p}(\R^N),\|\cdot\|^p,D)$ with the group of shifts $D=\{u\mapsto u(\cdot+y)\}_{y\in\R^N}\}$. 
\item $(\mathcal D^{1,p}(G),\|\sqrt{L(u)}\|_{p}^p,D')$, where $G$ is a Carnot group of homogeneous dimension $Q$ with invariant subelliptic Lagrangean $L(u)=\sum_i|du(X_i)|$, $X_i$ are generators of the correspondent Lie algebra, $1<p<Q$, and $D'$ is a product group of the actions of left group shifts and of the group of dilation actions $u\mapsto t^\frac{Q-p}{p}u\circ\delta_t$, $t>0$, where  $\delta_t: G\to G$, $t\in(0,\infty)$ are homogeneous dilations on $G$.
In particular, $(\mathcal D^{1,p}(\R^N),\|\nabla \cdot\|_p^p, D')$,
$1<p<N$, where $D'$ is a product group of Euclidean shifts and of the group of dilation actions $u\mapsto t^\frac{N-p}{p}u(t\cdot)$, is a dislocation space.
 \end{enumerate}
\begin{definition} 
Let $X$ be a Banach space and let $D$ be a group of linear  isometries on $X$. 
One says that a sequence $u_k\in X$ {\em converges $D$-weakly} to $u\in X$ (to be denoted $u_k\cw u$) if
for any sequence $g_k\in D$, $g_k(u_k-u)\rightharpoonup 0$.
\end{definition}

\begin{lemma}
\label{weak}
Let $(X,F,D)$ be a dislocation space. If $F(u_k)\to 0$ then $u_k\cw 0$.
\end{lemma}
\begin{proof}
Since $F(u_k)\le 1$ for all $k$ sufficiently large, $\|u_k\|\le 1$.
On every weakly convergent renumbered subsequence of $u_k$,
$F(\wlim u_k)\le \lim F(u_k)=0$ and consequently $u_k\rightharpoonup 0$. 
Since $D$ preserves $F$, the same conclusion applies to $g_ku$ for any sequence $g_k\in D$.
\end{proof}

\begin{theorem}
\label{abstractcc}
Let $(X,F,D)$ be a dislocation space. If $u_{k}\in X$ is a bounded sequence, then there
exists a set $\N_0\subset\N$, $w^{(n)}  \in X$, sequences $\{g_{k} ^{(n)}\}_{k \in\N} \subset D$ with $g_{k} ^{(1)} =id$, satisfying \eqref{separates}, $n\in\N_0$,
such that for a renumbered subsequence,
\begin{eqnarray}
\label{w_n} &&w^{(n)}=\wlim {g_{k} ^{(n)}}^{-1}u_k,
\\
\label{norms} &&\sum_{n \in \N_0} F(w ^{(n)}) \le
\limsup F(u_k)
\\
\label{BBasymptotics} &&u_{k}  -  \sum_{n\in\N_0}  g_{k} ^{(n)}
w^{(n)}  \cw 0,
\end{eqnarray}
where the series $\sum_{n\in\N_0}  g_{k} ^{(n)} w^{(n)}$ converges
uniformly in $k$ in the sense that
\begin{equation}
\label{uniformly}
\sup_{k\in\N}F\left(\sum_{n\ge m}  g_{k} ^{(n)} w^{(n)}\right)\to 0\mbox{ as } m\to\infty.
\end{equation}
\end{theorem}
\noindent
\begin{proof} 
\par\noindent 1. Once \eqref{separates} is proved, and $w^{(n)}$ satisfy \eqref{w_n},  inequality \eqref{wBL} from Definition~\ref{DS} holds for every $M\in\N$ and thus the series in (\ref{norms}) converges.
\par\noindent
2. Observe that if $u_k\cw 0$, the theorem is verified with
$\N_0=\emptyset$. Otherwise consider the expressions of
the form
\begin{equation}
\label{def_w1}
w^{(1)}=:\wlim
{g_k^{(1)}}^{-1}u_k.
\end{equation}
The sequence $u_k$ is bounded, $D$ is a set of isometries, so  
the sequence in (\ref{def_w1}) is bounded and
thus, for any choice of $g_k\in D$, it has a weakly convergent
subsequence. Since we assume that $u_k$ does not converge
$D$-weakly to zero, there exists necessarily a renumbered sequence
$g_k^{(1)}$ that yields a non-zero limit in (\ref{def_w1}).
\par
Let
\begin{equation*}
%\label{def_v1}
v_k^{(1)}= u_k - g_k^{(1)} w^{(1)},
\end{equation*}
and observe by (\ref{def_w1}) that
\begin{equation}
\label{g1v1} 
{g_k^{(1)}}^{-1}v_k^{(1)}={g_k^{(1)}}^{-1}(u_k-
w^{(1)})\rightharpoonup 0.
\end{equation}
If $v_k^{(1)} \cw 0$, the theorem is verified with $\N_0=\{1\}$.
If not -- we repeat the argument above -- there exist,
necessarily, a sequence $g_k^{(2)}\in D$ and a $w^{(2)}\neq 0$
such that, on a renumbered subsequence,
\begin{equation*}
%\label{def_w2}
{g_k^{(2)}}^{-1}
v_k^{(1)} \rightharpoonup w^{(2)}.
\end{equation*}
Let us set
\begin{equation*}
%\label{def_v2}
v_k^{(2)}= v_k^{(1)} - g_k^{(2)} w^{(2)}.
\end{equation*}
Then we will have an obvious analog of (\ref{g1v1}):
\begin{equation}
\label{w2-1} 
 {g_k^{(2)}}^{-1} v_k^{(2)}=
 {g_k^{(2)}}^{-1}(v_k^{(1)}  -
 w^{(2)} )\rightharpoonup 0.
\end{equation}
If we assume that
\begin{equation*}
%\label{w2-2}
{g_k^{(1)}}^{-1}g_k^{(2)}\not\rightharpoonup 0,
\end{equation*}
then by (\ref{w2-1}) and \eqref{newii},
\begin{equation*}
%\label{w2-3}
{g_k^{(1)}}^{-1} (v_k^{(1)} - g_k^{(2)}w^{(2)})\rightharpoonup 0,
\end{equation*}
which, due to (\ref{g1v1}), yields
\begin{equation}
\label{w2-4}
{g_k^{(1)}}^{-1} g_k^{(2)}w^{(2)}\rightharpoonup 0.
\end{equation}
We now use (\ref{newii}) again to replace in (\ref{w2-4})
$ {g_k^{(1)}}^{-1}$ with
${g_k^{(2)}}^{-1}$, which results in
\begin{equation}
 w^{(2)}\rightharpoonup 0,
\end{equation}
 which cannot be true since we
assumed $w^{(2)}\neq 0$. From this contradiction follows
\begin{equation}
\label{separation_1_02}
{g_k^{(1)}}^{-1} g_k^{(2)} \rightharpoonup 0.
\end{equation}
Then we also have
\begin{equation*}
%\label{separation_2_1}
{g_k^{(2)}}^{-1} g_k^{(1)} \rightharpoonup 0.
\end{equation*}
Indeed, it this were false, then from \eqref{newii} and \eqref{separation_1_02} we
have on a subsequence
\begin{equation*}
%\label{separation_2_1}
id = {g_k^{(2)}}^{-1} g_k^{(1)} {g_k^{(1)}}^{-1} g_k^{(2)} \rightharpoonup 0,
\end{equation*}
which is obviously false.
Recursively we define:
\begin{equation}
\label{def:next_v} v_k^{(n)}\eqdef v_k^{(n-1)}-g_k^{(n)} w^{(n)} =
u_k - g_k^{(1)} w^{(1)} - \dots - g_k^{(n)} w^{(n)},
\end{equation}
where
\begin{equation*}
%\label{def_wn}
w^{(n)}=\wlim
{g_k^{(n)}}^{-1} v_k^{(n-1)},
\end{equation*}
calculated on a successively renumbered subsequence. We
subordinate the choice of $g_k^{(n)}$ and thus extraction of this
subsequence for every given $n$ to the following requirements. For
every $ n\in\N$ we set
\begin{equation*}
%\label{setW}
W_n=\{w\in H\setminus\{0\}:\;\exists g_j\in D,
\{k_j\}\subset\N: g_j^{-1}v_{k_j}^{(n)}\rightharpoonup
w\},
\end{equation*}
and
\begin{equation*}
%\label{t_n}
t_n=\sup_{w\in W_n}F(w).
\end{equation*}
Note that $t_n<\infty$, since all operators involved at all steps
leading to the definition of $W_n$ have uniform bounds. 

If for
some $n$, $t_n=0$, the theorem is proved with $\N_0=\{1,\dots,n-1\}$. Otherwise, we choose a
$w^{(n+1)}\in W_n$ such that
\begin{equation}
\label{w-t} F(w^{(n+1)})\ge\frac12 t_n
\end{equation}
and the sequence $g_k^{(n+1)}$ is chosen so that on a subsequence
that we renumber,
\begin{equation}
\label{next_w}
{g_k^{(n+1)}}^{-1}v_k^{(n)}
\rightharpoonup w^{(n+1)}.
\end{equation}
An argument analogous to the one brought above for $n=1$ shows
that
\begin{equation}
\label{separation_p_q} {g_k^{(p)}}^{-1}g_k^{(q)}\rightharpoonup 0 \mbox{
whenever } p \neq q, p,q\leq n.
\end{equation}
This allows to deduce immediately (\ref{w_n}) from (\ref{next_w}) as well as \eqref{norms}.
%Note that the sequence $t_n$ is non-increasing.
%
From (\ref{wBL}) and(\ref{w-t}) follows
\begin{equation*}
%\label{SNE5}
\sum_{n\ge 2}t_n\le 2F(u_k).
\end{equation*}
Let $\varphi_i$, $i\in\N$, be a normalized basis in $X^*$. Then by definition
of $W_n$,
\begin{equation*}
\limsup_k \sum_i2^{-i} \sup_{g\in D}
\langle gv^{(n)}_k,\varphi_i\rangle^2\le 4t_n^2,\;n\in\N.
\end{equation*}
Let $k(n)$ be such that
\begin{equation}
\label{k_n} \sum_i2^{-i} \sup_{g\in D}
\langle gv^{(n)}_{k(n)},g\varphi_i\rangle^2\le 8t_n^2, \;n\in\N.
\end{equation}
This implies that
\begin{equation*}
\sup_{g\in D} \langle gv^{(n)}_{k(n)},\varphi\rangle\to 0
\end{equation*}
for any $\varphi$ that is a linear combination of $\varphi_i$,
and an elementary density argument extends this relation to any
$\varphi\in X^*$, so that
\begin{equation*}
v^{(n)}_{k(n)}\cw 0
\end{equation*}
as $n\to\infty$. Instead of $k(n)$ selected for each $n$ from the
index set of a renumbered subsequence of $u_k$ (that was produced
by successive extractions), we will now use the correspondent
index (preserving the notation $k(n)$) from the original
enumeration of $u_k$. (This change of enumeration affects also the
terms $g^{(j)}_{k(n)}$, $j=1,\dots,n$, in the definition
(\ref{def:next_v}) of $v^{(n)}_{k(n)}$.) Then we conclude that
\begin{equation*}
v^{(n)}_{k(n)}=u_{k(n)}-\sum_{j\le n} g^{(j)}_{k(n)}w^{(j)}\cw 0.
\end{equation*}
Since
the final extraction is a subsequence of the sequence in
(\ref{separation_p_q}), (\ref{separates}) follows.

Note that $k(n)$ can be chosen in (\ref{k_n}) arbitrarily large,
and in particular large enough so that the series $\sum_j
g^{(j)}_{k(n)}w^{(j)}$ is uniformly convergent in the sense of \eqref{uniformly} due to
(\ref{norms}) and \eqref{separates}, and therefore (\ref{BBasymptotics}) follows.
Indeed, one can always choose a subsequence of $g^{(m+1)}_{k}$ so that, by \eqref{w1BL}, 
\begin{eqnarray*}
\left|F\left(\sum_{n=1}^{m+1}
g^{(n)}_{k}w^{(n)}\right)-F\left(\sum_{n=1}^{m}
g^{(n)}_{k}w^{(n)}\right)-F(w^{(m+1)})\right|\le 2^{-k-m}.
\end{eqnarray*}
Finally, if $w^{(1)}=\wlim u_k\neq 0$, we could have chosen
$g_k^{(1)}=id$ at the first step. If $\wlim u_k=0$, we renumber
the terms in expansion by $n=2,3,\dots$ and set $g_k^{(1)}=id$,
$w^{(1)}=0$.
\end{proof}

\section{Cocompactness and minimizers}
In this section we give a functional-analytic formalization of the minimization reasoning of P.-L.Lions (\cite{PLL1a}) in cocompactly imbedded dislocation spaces.

\begin{definition} One says that a continuous imbedding of a Banach space $X$ into 
a Banach space $Y$ is cocompact relatively to a  group  $D$ of isometric linear operators on $X$ if every $D$-weakly convergent sequence $u_k\in X$ converges in $Y$.
\end{definition}
Note that it does not follow from this definition that the quotient $X/D$ is compactly imbedded into $Y$. 
If $D=\{id\}$, the cocompact imbedding becomes compact. 
\par\noindent
The following imbeddings are cocompact. 
\begin{enumerate}
\item $W^{1,p}(\R^N)$, is cocompactly imbedded into $L^q(\R^N)$,
relatively to the group of lattice shifts $u\mapsto u(\cdot+y)$, $y\in\Z^N$, when $p<q<\frac{Np}{N-p}$ for $N>p$ or $q>p$ for $N\le p$. 
\item Let $M$ be a complete $N$-dimensional Riemannian manifold cocompact with respect to a subgroup $G$ of its isometry group
$Iso(M)$, that is, there exists a compact set $V\subset M$ such that $\cup_{\eta\in G}\eta V=M$. Then,  $W^{1,p}(M)$ with the invariant norm $\|u\|^p=\int(|du|^p+|u|^p)d\mu$ is cocompactly imbedded into $L^p(M)$ for the same values of $p$ as above, relatively to the group 
$\{u\mapsto u\circ\eta\}_{\eta\in G}$.
\item Let $G$ be a Carnot group of homogeneous dimension $Q$. Then 
$\mathcal D^{1,p}(G)$, $p<Q$, is cocompactly imbedded into $L^{p^*}(G)$
where $p^*=\frac{pQ}{Q-p}$, relatively to a product group of left shifts and 
discrete dilation action $u\mapsto 2^{\frac{Q-p}{p}j}u\circ\delta_{2^j}$, $j\in\Z$. In particular, $\mathcal D^{1,p}(\R^N)$ is cocompactly imbedded into $L^{p^*}(\R^N)$ for $p<N$.
\end{enumerate}
The Euclidean case of the statements above, with the group of $\R^N$-shifts, and, in the limit Sobolev case, with the continuous dilation group, is due to Lieb \cite{Lieb} and Lions (\cite{PLL1a,PLL2a}). The proof in the case of a manifold and of discrete dilations is found in \cite{ccbook} for $p=2$. 
The general case can be proved by direct analogy with those.

In what follows assume that 
\vskip2mm\noindent
({\bf A}) $(X,F,D)$ is a dislocation space, $(Y,G)$ is as in Lemma~\ref{prop:OS}, $G\in C(X)$,
$X$ is continuously imbedded into $Y$, and $D$ has an extension into $Y$ such that $G\circ g=G$ for all $g\in D$. Moreover, $G$ satisfies the Brezis-Lieb property:
\begin{equation}
\label{BL}
G(u_k)-G(u)-G(u-u_k)\to 0 \text{ whenever } u_k\rightharpoonup u \text{ in } X.
\end{equation}

\begin{lemma}
\label{lem:cont} Let $(Y,G,D)$ satisfy assumption ({\bf A}).
Then the function $G$ is continuous in $Y$.
\end{lemma}
\begin{proof} Let $u_k\to u$ in $Y$.
By Lemma~\ref{prop:OS}, $G(u_k-u)\to 0$. By \eqref{BL}, 
$$\lim G(u_k)-G(u)=\lim G(u_k-u)=0.$$
\end{proof}

Let
\begin{equation}
\label{isop}
c_t\eqdef\inf_{u\in X: G(u)=t} F(u), \qquad t>0.
\end{equation}

\begin{proposition} 
\label{prop:wsubadd}
Assume ({\bf A}) and 
\par\noindent({\bf B}) $D$ contains a subsequence $g_k\rightharpoonup 0$.
\par\noindent
Then for any $\tau\in[0,t]$,
\begin{equation}
\label{wsubadd}
c_t\le c_\tau+c_{t-\tau}.
\end{equation}
\end{proposition}
\begin{proof}
Fix $\epsilon>0$ and let 
$v,w\in X$ satisfy respectively
$G(v)=\tau$, $F(v)\le c_\tau+\epsilon/2$
and
$G(w)=t-\tau$, $F(w)\le c_{t-\tau}+\epsilon/2$.
Let $g_k\rightharpoonup 0$ and let $u_k=v+g_kw$.
Then by \eqref{BL}, $G(u_k)\to G(v)+G(w)=t$. On the other hand, by \eqref{w1BL} 
$F(u_k)\to F(v)+F(w)\le c_\tau+c_{t-\tau}+\epsilon$.
Since $\epsilon$ is arbitrary, this implies \eqref{wsubadd}.
\end{proof}
We show existence of constrained minima under assumptions of the strict inequality in \eqref{wsubadd} and of cocompactness.
This is a functional-analytic formalization of the analogous results of Lions.

\begin{lemma}
\label{lem:coercivity}Assume ({\bf A}) and let the embedding of $X$ into $Y$ be cocompact. 
Then for all $a,b>0$
\begin{equation}
\inf_{G(au)>b}F(u)>0.
\label{coercivity}
\end{equation}
\end{lemma}
\begin{proof}
Assume that there is a sequence $u_k\in X$ such that $F(u_k)\to 0$, while
$G(au_k)>b$. By Lemma~\ref{weak}, $au_k\cw 0$ and by cocompactness of imbedding $au_k\to 0$ in $Y$. By Lemma~\ref{lem:cont}, $G$ is continuous and therefore
$G(au_k)\to 0$, a contradiction. 
\end{proof}

\begin{theorem} 
\label{Lions} Assume {\rm({\bf A})} and  
({\bf B}).
Then for every minimizing sequence $u_k$ for \eqref{isop}, $t>0$ there exists a sequence $g_k\in D$ such that $g_ku_k$ converges $D$-weakly to a point of minimum. 
if and only if for every $\tau\in(0,t)$
\begin{equation}
\label{subadd}
c_t<c_{\tau}+c_{t-\tau}.
\end{equation}
\end{theorem}

\begin{proof} Note that by  \eqref{coercivity} we have $c_t>0$.
\par\noindent {\it Sufficiency.} Assume \eqref{subadd}.
Let $u_k\in X$ be a minimizing sequence, that is, $F(u_k)=c_t$
and $G(u_k)\to t$.
If $u_k\cw 0$, then $u_k\to 0$ in $Y$ by cocompactness, and $G(u_k)\to 0$ by Lemma~\ref{lem:cont}, which contradicts $c_t>0$.
Consequently, there exists a sequence $g_k\in D$ such that, on a renamed subsequence, $g_ku_k\rightharpoonup w^{(1)}\neq 0$ in $X$. By invariance of the problem, $g_ku_k$ is also a minimizing sequence that we now rename as $u_k$. Let $g_k^{(n)}$, $w^{(n)}$ be as provided by Theorem~\ref{abstractcc}. 
By \eqref{norms} 
$$
\sum_nF(w^{(n)})\le c(t),
$$
and from iteration of \eqref{BL} and cocompactness of the imbedding follows
$$
\sum_nG(w^{(n)})=t.
$$
Let $G(w^{(n)})=\tau_n$. Then
$$
\sum_nc_{\tau_n}\le c_t,
$$
which by \eqref{subadd} is false unless all but one of the values $\tau_n$ is zero.
Since $\tau_1\neq 0$, we conclude that
$u_k-w^{(1)}\cw 0$. By cocompactness $u_k\to w^{(1)}$ in $Y$ and by continuity of $G$, $G(w^{(1)}=t$. By weak lower semicontinuity, $F(u)\le c_t$. Since $c_t$ is the infimum over functions with $G(u)=t$, $w^{(1)}$ is necessarily a minimizer.

\noindent{\em Necessity.} Assume now that \eqref{subadd} does not hold for some $0<\tau<t$.
By \eqref{wsubadd} this implies $c_\tau+c_{t-\tau}=c_t$.
Let  
$v_n,w_n\in X$ satisfy respectively
$G(v_n)=\tau$, $F(v_n)\le c_\tau+1/n$
and
$G(w_n)=t-\tau$, $F(w_n)\le c_{t-\tau}+1/n$, $n\in\N$.
Let $g_k\rightharpoonup 0$ and let $u_{nk}=v_n+g_kw_n$.
Then for every $n$ there exists $k_n$ such that for all 
$k\ge k_n$,  
$\sup_{k\ge k_n}|G(u_{nk})-t|\to 0$ by \eqref{BL}, and 
$\sup_{k\ge k_n}|F(u_{nk})-c_{t}|\to 0$ by \eqref{w1BL}, as $n\to\infty$.
Without loss of generality, $v_n\rightharpoonup v\neq 0$ and
$w_n\rightharpoonup w\neq 0$ (if one of $v_n$ and $w_n$ is $D$-weakly convergent to zero, then $\tau=0$ or $\tau=t$).
Let now $\psi_j$, $j\in\N$ be a basis in $X^*$. Then
$$
\sum_{j\in\N}\frac{|\langle\psi_j,g_{k'_n}w_n\rangle|^2}{2^j}\to 0
$$
if only one chooses $k'_n\ge k_n$ sufficiently large. This implies that 
$g_{k'_n}w_n\rightharpoonup 0$.
A similar argument allows to select a further subsequence
such that $g_{k''_n}^{-1}v_n\rightharpoonup 0$.
Consequently, $\wlim(v_n+g_{k''_n}w_n)=v\neq 0$, 
while $\wlim(g_{k''_n}^{-1}v_n+w_n)=w\neq 0$. Thus we have constructed a minimization sequence that is not $D$-weakly convergent.
\end{proof}
Note that that the proof of sufficiency does not require condition ({\bf B}).

\begin{theorem}
\label{penalty}
Let $(X,F,D)$ be a dislocation space, and assume \eqref{subadd}, ({\bf A})  and ({\bf B}).
Let $f,g:X\to\R$ be nonnegative weakly continuous functions, at least one of them is positive for $u\neq 0$, and let
\begin{equation}
\label{isop2}
c'_t\eqdef\inf_{G(u)+g(u)=t}(F(u)-f(u)), t>0.
\end{equation}
If for every $\tau\in(0,t)$
 
\begin{equation}
\label{subadd2}
c'_t<c'_\tau+c_{t-\tau},
\end{equation}
then every minimizing sequence for \eqref{isop2} converges $D$-weakly to a point of minimum.
\end{theorem}
Using an argument repetitive of that in Proposition~\ref{prop:wsubadd} one can easily see that $c'_t\le c'_\tau+c_{t-\tau}$ for any $\tau\in(0,t)$, $t>0$, so the role of condition \eqref{subadd2} is similar to that of \eqref{subadd}.
\begin{proof}
The proof is analogous to the proof of Theorem~\ref{Lions} and to the similar statement in \cite{PLL1a}.
Let $u_k\in X$ be a minimizing sequence, that is, $(F-f)(u_k)=c'_t$
and $(G+g)(u_k)\to t$.
Let $g_k^{(n)}$, $w^{(n)}$ be as provided by Theorem~\ref{abstractcc}. 
By iteration of \eqref{BL}, taking into account cocompactness,
$$
g(w^{(1)})+\sum_nG(w^{(n)})=t.
$$
Let $G(w^{(1)})+g(w^{(1)})\eqdef\tau_1$ and $G(w^{(n)})\eqdef\tau_n$, $n\ge 2$,
so that 
$\sum \tau_n\le t$.
From 
By \eqref{norms} one has also
$$
\sum_nF(w^{(n)})-f(w^{(1)})\le c'(t),
$$
which implies 
$$
c'_{\tau_1}+\sum_{n\ge 2}c_{\tau_n}\le c'_t.
$$
This contradicts \eqref{subadd} and \eqref{subadd}
unless all but one of the values $\tau_n$ is zero. Assume that $\tau_m=1$ for some $m\ge 2$.
Then $c_t\le c'_t$, which is false (the opposite strict inequality follows from the from substitution of the minimizer of \eqref{isop} into \eqref{isop2}). Consequently, $u_k-w^{(1)}\cw 0$, $(G+g)(w^{(1)})=t$ and $(F-f)(w^{(1)})\le c'_t$, so $w^{(1)}$ is necessarily a minimizer.
\end{proof}
%Consider now unconstrained minimization.
%Let $\Psi\eqdef F-G$, where $X,F$ and $Y,G$ are as in %Lemma~\ref{prop:OS}. 
%Note that 
%$$
%\inf_X\Psi=\inf_{t>0} (c_t-t),
%$$ 
%where $c_t$ is given by \eqref{isop}.
%Note that since $F(0)=G(0)=0$, the infimum is in general nonnegative. 

%We have immediately
%\begin{proposition}
%\label{unconstrained}
%Assume conditions of Theorem~\ref{Lions}. 
%If $\inf_{t>0} (c_t-t)$ is attained at some $t_0>0$, then $\inf_X\Psi$
%is attained at some $u\in X$ with $F(u)=c_{t_0}$, $G(u)=t_0$.
%\end{proposition}
%Assumptions of Proposition~\ref{unconstrained} imply %$\inf_X\Psi>-\infty$, and are usually verified when $\inf_{t>0} (c_t-t)<0$.

\section{Flask subspaces}
Theorem~\ref{abstractcc} can be extended to certain subspaces of 
a dislocation space which are not $D$-invariant. 
\begin{definition} 
\label{def:flask}
Let $(X,F,D)$ be a dislocation space. A subspace $X_0$ of $X$ is called a flask subspace if $g_ku_k\rightharpoonup u$, $g_k\in D$, $u_k\in X_0$,
implies that $gu\in X_0$ for some $g\in D$.
\end{definition}
\begin{proposition}
\label{prop:flask} let $(X,F,D)$ be a dislocation space with a flask subspace $X_0$ and let $u_k\in X_0$ be a bounded sequence. Then Theorem~\ref{abstractcc} holds with $w^{(n)}\in X_0$.
\end{proposition}
\begin{proof}
Since $X_0$ is a flask subspace, $g_nw^{(n)}\in X_0$ for some $g_n\in D$, $n\in\N$. Set $\tilde w^{(n)}\eqdef g_nw^{(n)}$ and $\tilde g_k^{(n)}\eqdef g_k^{(n)}g_n^{-1}$, so that \eqref{BBasymptotics} holds 
with $\tilde w^{(n)}$ and $\tilde g_k^{(n)}$. It is easy to see that sequences $\tilde g_k^{(n)}$ satisfy \eqref{separates}. 
\end{proof}
\begin{corollary}
Proposition~\ref{prop:wsubadd}, Lemma~\ref{lem:coercivity}, 
Theorem~\ref{Lions} and Theorem~\ref{penalty} remain valid if $X$ in the statements \eqref{isop},\eqref{isop2} is replaced by a flask subspace $X_0$.
\end{corollary}

Flask subspaces are a functional-analytic generalization of the case $H^1_0(\Omega)$ with a flask domain $\Omega\subset\R^N$ in the sense of del Pino - Felmer \cite{PinoFelmer}. 

\begin{proposition}
\label{flask-suff}
Let $M$ be a complete Riemannian manifold, cocompact with respect to a subgroup $G$ of its isometry group $Iso(M)$. Let $\Omega\subset M$ be an open set with a piecewise smooth boundary.
If for every sequence $\eta_k\in G$ there exists $\eta\in Iso(M)$ such that
\begin{equation}
\label{eq:flask}
\liminf \eta_k(\Omega)\subset \eta(\Omega),
\end{equation}
then $W^{1,p}_0(\Omega)$, $p>1$, is a flask subspace of $W^{1,p}(M)$ relatively to the group $\{u\mapsto u\circ\eta\}_{\eta\in G}$.
\end{proposition}
\begin{proof}
First observe that, for arbitrary functions, if $u_k(x)\to u(x)$ and $u(x)\neq 0$, then necessarily $u_k(x)\neq 0$ for all $k$ sufficiently large. In other words, 
$$
\{u\neq 0\}\subset \liminf\{u_k\neq 0\}
$$
If $u_k\circ\eta_k\rightharpoonup u$ in $W^{1,p}(M)$, then $u_k\circ\eta_k$ converges almost everywhere as well, and we conclude from \eqref{eq:flask} that for some $\eta\in Iso(M)$, $u=0$ a.e. on $M\setminus \eta(\Omega)$. In order to apply Hedberg's trace theorem \cite{AdamsHedberg} (to regularized $u$) it remains to note that $u=0$ on $M\setminus (\overline{\eta(\Omega)}$ and, since $\partial\Omega$ is sufficiently smooth, $u=0$ on $\eta(\partial\Omega)$ as well, which yields $u\in W^{1,p}_0(\eta(\Omega))$. 
\end{proof}

\section{Compact imbeddings}
This section deals with abstract analogs of sufficient conditions for compactness of Sobolev imbeddings on unbounded domains (see e.g \cite{Adams70},\cite{CClark66}). 

\begin{proposition} 
\label{null}
Let $(X,F,D)$ be a dislocation space cocompactly imbedded into a Banach space $Y$. Assume (B). Let $X_0$ be a subspace of $X$.
If for every sequence $u_k\in X_0$ 
\begin{equation}
\label{novanish}
\{g_k\}\subset D, g_k\rightharpoonup 0\Rightarrow g_ku_k\rightharpoonup 0,
\end{equation}
then the imbedding of $X_0$  into $Y$ is compact.
\end{proposition}
\begin{proof}
By (B), a sequence $g_k\rightharpoonup 0$ exists. Then, since the sequence $g_ku_k$ is bounded, and $g_k$ are isometries, $u_k$ is a bounded sequence. Without loss of generality it suffices to assume that $u_k$ has the form \eqref{BBasymptotics}. Then condition 
\eqref{novanish} implies $u_k-\wlim u_k\cw 0$. Since the imbedding of $X$ into $Y$ is cocompact, this implies $u_k-\wlim u_k\to 0$ in $Y$.
\end{proof}

\begin{corollary} Let $M$ be a sub-Riemannian manifold of homogeneous dimension $Q$ cocompact with respect to $Iso(M)$. 
If $\Omega\subset M$ is an open set and if
for any sequence $\eta_k\in Iso(M)$ such that, for some $x_0\in M$, $\eta_k(x_0)$ has no convergent subsequence,
$$
\liminf \eta_k(\Omega) \text{ has measure zero},
$$ 
then $W_0^{1,p}(\Omega)$ is compactly imbedded into $L^q(\Omega)$, 
$1<p<q<p^*$.
\end{corollary}
\begin{proof}
Since the imbedding in question is cocompact, the statement follows from  Proposition~\ref{null} once we observe that the operator sequence $u\mapsto u\circ\eta_k$, with $\eta_k$ as above, is weakly convergent to zero.
Indeed, if it does not, then, necessarily, there exists a compact set $V\subset M$ such that, for a renamed subsequence, $\cup_k\eta_kV$ is a bounded set. Since $\eta_k$ are isometries, this yields, by Arzela-Ascoli theorem, that a subsequence of $\eta_k$ is convergent uniformly on compact sets, and in particular, $\eta_k(x_0)$ is convergent, a contradiction.
\end{proof}
The following statement generalizes the well known compactness for subspaces of radial functions (e.g. \cite{EstebanLions83}).
\begin{theorem}
Let $(X,F,D)$ be a dislocation space cocompactly imbedded into a Banach space $Y$. Let $C$ be a group of linear automorphisms of $X$ that preserves $F$, such that for every $c\in C\setminus\{id\}$ and 
for every sequence $g_k\in D$,
$g_k\rightharpoonup 0$, 
\begin{equation}
\label{everyc}
g_k^{-1}cg_k\rightharpoonup 0. 
\end{equation}
Let
$$
X_C\eqdef\{u\in X: cu=u, c\in C\}.
$$
Then the imbedding of the subspace $X_C$ into $Y$ is compact.
\end{theorem}
\begin{proof}
Let $u_k$ be a bounded sequence in $X_C$ and consider its expansion \eqref{BBasymptotics}. Then for any $c\in C$, $c^{-1}u_k=u_k$ and therefore
\begin{equation}
u_k-\sum_n cg_k^{(n)}w^{(n)}
\end{equation}
Assume that there is at least one term $w^{(n)}\neq 0$, with $n>2$, say, with $n=2$.
Then by \eqref{everyc},  
$$
{g_k^{(2)}}^{-1}cg_k^{(2)}\rightharpoonup 0, c\in C\setminus\{id\},
$$
for every $c,c'\in C$, $c'\neq c$ 
$$
(c'g_k^{(2)})^{-1}cg_k^{(2)}\rightharpoonup 0,
$$
and, furthermore
$$
(cg_k^{(2)})^{-1}u_k\rightharpoonup w^{(2)}, c\in C.
$$
Let $M\in\N$ and let $C_M$ be any subset of $C$ with $M$ elements. Then, by \eqref{wBL},
$$F(u_k)\ge \sum_{c\in C_M}F(w^{(2)})=MF(w^{(2)}).$$
 Since $M$ is arbitrary and the left hand side is bounded, we arrive at a contradiction. Consequently, $u_k\cw w^{(1)}$. Since the imbedding of $X_C$ into $Y$ is cocompact, $u_k\to w^{(1)}$ in $Y$.
\end{proof}

\begin{center}  
{\bf Acknowledgments.}
\end{center}
The topic of this paper was suggested, at different times, by Moshe Marcus and Adimurthi. The author thanks Vladimir Maz'ya for a discussion concerning the use of Hedberg trace theorem in the context of Proposition~\ref{flask-suff}. The paper was written as the author was visiting Technion -- Haifa Institute of Technology and he thanks the mathematics faculty there for their warm hospitality.


\begin{thebibliography}{99}
\bibitem{Adams70}
Adams, R.A., Compact imbedding theorems for quasibounded domains, 
Trans.Amer.Math.Soc., {\bf 148}, 445-459 (1970)

\bibitem{AdamsHedberg} Adams, D.R.; Hedberg, L.I.,
Function Spaces and Potential Theory. Springer-Verlag, 1995.
\bibitem{BL} Brezis, H.; Lieb, E.,
A relation between pointwise convergence of functions and
convergence of functionals. {\it Proc. Amer. Math. Soc.} {\bf 88}
(1983), 486--490.
%
%
\bibitem{BrezisCoron} Brezis, H; J.M. Coron;
Convergence of solutions of H-systems or how to blow bubbles, {\it
Archive Rat. Mech. Anal.} {\bf 89} (1985), 21--56. 
%
\bibitem{CClark66} Clark C., An embedding theorem for function spaces, 
Pacific J.Math. {\bf 19}, 243-251 (1966)
\bibitem{EstebanLions83}
Esteban M., Lions P.-L., 
A compactness lemma, Nonlinear Anal. {\bf 7}, 381-385 (1983).
\bibitem{Lieb} Lieb, E., On the lowest eigenvalue of the Laplacian for the intersection of two domains.
Invent. Math. {\bf 74}, 441-448 (1983)
\bibitem{PLL1a} Lions P.-L., The concentration-compactness principle in the 
calculus of variations. The locally compact case, part 1. 
Ann.Inst.H.Poincare, Analyse non lin\'eaire {\bf 1}, 109-1453 (1984)
\bibitem{PLL1b} Lions P.-L., The concentration-compactness principle in the 
calculus of variations. The locally compact case, part 2. 
Ann.Inst.H.Poincare, Analyse non lin\'eaire {\bf 1}, 223-283 (1984)
\bibitem{PLL2a} Lions P.-L., The concentration-compactness principle in the 
calculus of variations. The Limit Case, Revista Matematica Iberoamericana,
Part 1, {\bf 1.1} 145-201 (1985)
\bibitem{PLL2b} Lions P.-L., The concentration-compactness principle in the 
calculus of variations. The Limit Case, Revista Matematica Iberoamericana,
Part 2, {\bf 1.2} 45-121 (1985)
\bibitem{Lions87} Lions, P.-L., Solutions of Hartree-Fock equations for Coulomb systems,
{\it Comm.~Math.~Phys.} {\bf 109}  (1987), 33--97.
\bibitem{Mazja1}Maz'ja, V.G., Sobolev Spaces, Springer-Verlag (1985)  
\bibitem{PinoFelmer} del Pino, M., Felmer, P., Least energy solutions for
elliptic equations in unbounded domains, Proc. Royal Soc. Edinburgh {\bf
126A}, 195-208 (1996)
\bibitem{Struwe84} Struwe, M., A global compactness result for elliptic boundary value problems involving limiting nonlinearities.
{\it Math. Z.} {\bf 187} (1984), 511--517.
\bibitem{ccbook} Tintarev K., Fieseler K.-H., Concentration compactness: functional-analytic grounds and applications, Imperial College Press 2007.
\end{thebibliography}
\end{document}